\documentclass[11pt,letter,english]{amsart}

\font \sevenrm=cmr7
\font \fiverm=cmr5

\usepackage{amsfonts}
\usepackage{amssymb}
\usepackage{amsmath}
\usepackage[all]{xy}  
\usepackage{color}
\usepackage{graphicx}
\usepackage{xspace}

 \newcommand{\nc}{\newcommand}
\usepackage{mathrsfs} 

 {\everymath{\displaystyle\everymath{}}\array}%
 {\endarray}

\newtheorem{thm}{Theorem}
\newtheorem{exam}{Example}

\newtheorem{cor}[thm]{Corollary}

\newtheorem{lem}[thm]{Lemma}
\newtheorem{prop}[thm]{Proposition}
\newtheorem{defn}{Definition}

\newtheorem{problem}[thm]{Probleme}

\nc{\comment}[1]{[[{\tt #1}]] }
\nc{\Cal}[1]{{\mathcal {#1}}}
\nc{\mop}[1]{\mathop{\hbox {\rm #1} }\nolimits}
\nc{\gmop}[1]{\mathop{\hbox {\bf #1} }\nolimits}
\nc{\smop}[1]{\mathop{\hbox {\sevenrm #1} }\nolimits}
\nc{\ssmop}[1]{\mathop{\hbox {\fiverm #1} }\nolimits}
\nc{\mopl}[1]{\mathop{\hbox {\rm #1} }\limits}
\nc{\smopl}[1]{\mathop{\hbox {\sevenrm #1} }\limits}
\nc{\ssmopl}[1]{\mathop{\hbox {\fiverm #1} }\limits}
\nc{\frakg}{{\frak g}}
\nc{\g}[1]{{\frak {#1}}}
\def \restr#1{\mathstrut_{\textstyle |}\raise-6pt\hbox{$\scriptstyle #1$}}
\def \srestr#1{\mathstrut_{\scriptstyle |}\hbox to
  -1.5pt{}\raise-4pt\hbox{$\scriptscriptstyle #1$}}
\nc{\wt}{\widetilde} \nc{\wh}{\widehat}
\nc{\redtext}[1]{\textcolor{red}{\tt [[#1]]}}
\nc{\bluetext}[1]{\textcolor{blue}{#1}}
\nc\fleche[1]{\mathop{\hbox to #1 mm{\rightarrowfill}}\limits}
\nc{\ignore}[1]{}
\def\semi{\mathrel{\times}\kern -.85pt\joinrel\mathrel{\raise
    1.4pt\hbox{${\scriptscriptstyle |}$}}}
\nc\R{{\mathbb R}}
\nc\N{{\mathbb N}}
\nc\inver{^{-1}}
\nc\point{\hbox{\bf .}}
\nc\un{\hbox{\bf 1}}
\def\link#1#2{\raise -2pt\hbox{$\scriptstyle #1-\!\!-\!\!- #2$}}
\def\slink#1#2{\raise -1.5pt\hbox{$\scriptscriptstyle #1-\!\!\!-\!\!\!- #2$}}

\def\bluetext#1{\textcolor[rgb]{0.00,0.00,1.00}{#1}}

\def\redtext#1{\textcolor[rgb]{0.98,0.00,0.00}{#1}}

\definecolor{Light}{gray}{0.85}

\def\abs#1{\left\vert #1 \right\vert}
\def\allpoly{\mbox{$\re\langle X \rangle$}}
\def\allpolyell{\mbox{$\re^{\ell}\langle X \rangle$}}
\def\allpolyx0degn{\mbox{$P_n$}}

\def\allseries{\mbox{$\re\langle\langle X \rangle\rangle$}}

\def\allseriesell{\mbox{$\re^{\ell} \langle\langle X \rangle\rangle$}}

\def\allseriesellGC{\mbox{$\re^{\ell}_{GC}\langle\langle X \rangle\rangle$}}

\def\bfem#1{{\bf \em #1}} 

\def\bull{\rule{0.08in}{0.08in}} 

\def\diag{{\rm diag}}

\def\eqref#1{(\ref{#1})} 

\def\expup{{\rm e}}

\def\Lpm{L_{\mathfrak{p}}^m}
\def\Lpe{L_{\mathfrak{p},e}}
\def\Lpme{L^m_{\mathfrak{p},e}}

\def\mbf#1{\hbox{\mathversion{bold}$#1$}} 

\def\norm#1{\left\Vert#1\right\Vert}

\def\re{{\mathbb R}} 

\def\sameau{\rule[0.017in]{0.2in}{0.012in}}

\def\begals{\[\begin{aligned}}
\def\endals{\end{aligned}\]}
\def\begce{\begin{center}}
\def\endce{\end{center}}
\def\begar{\begin{array}}
\def\endar{\end{array}}
\def\begeq{\begin{equation}}
\def\endeq{\end{equation}}
\def\begdi{\begin{displaymath}}
\def\enddi{\end{displaymath}}
\def\begdis{\begin{eqnarray*}}
\def\enddis{\end{eqnarray*}}
\def\begeqa{\begin{eqnarray}}
\def\endeqa{\end{eqnarray}}
\def\begdes{\begin{description}}
\def\enddes{\end{description}}
\def\begit{\begin{itemize}}
\def\endit{\end{itemize}}
\def\begen{\begin{enumerate}}
\def\enden{\end{enumerate}}
\def\beglar{\left[\begin{array}}
\def\endrar{\end{array}\right]}
\def\begle{\begin{lem}}
\def\endle{\end{lem}}
\def\begde{\begin{defn}}
\def\endde{\end{defn}}
\def\begth{\begin{thm}}
\def\endth{\end{thm}}
\def\begco{\begin{cor}}
\def\endco{\end{cor}}
\def\begprop{\begin{prop}}
\def\endprop{\end{prop}}
\def\begex{\begin{exam}}
\def\endex{\end{exam}}

\def\begres{\noindent{\bf Remarks}:\begin{enumerate}}
\def\endres{\end{enumerate} \par}
\def\begpr{\noindent{\em Proof:}$\;\;$}
\def\endpr{\hfill\bull \vspace*{0.15in}}
\def\begtab{\begin{tabular}}
\def\endtab{\end{tabular}}
\def\rref#1{(\ref{#1})}
\def\lbracedef{\left\{\begin{array}{@{\hspace*{2pt}}c@{\hspace*{3pt}}c@{\hspace*{3pt}}l}} 

\def\allseriesA{\mbox{$\re\!\ll\!A\!\gg$}}
\def\allseriesA'{\mbox{$\re\!\ll\!A'\!\gg$}}

\begin{document}

\title[Evaluating Generating Functions for Periodic Multiple Polylogarithms]
{\vspace*{-0.7in}Evaluating Generating Functions for Periodic Multiple Polylogarithms via Rational Chen--Fliess Series}

\author{Kurusch Ebrahimi-Fard}
\address{ICMAT, C/ Nicol\'{a}s Cabrera, no.~13-15, 28049 Madrid, Spain {\tiny{On leave from UHA, Mulhouse, France.}}}
	\email{kurusch@icmat.es}
	\urladdr{http://www.icmat.es/kurusch/personal}
	
\author{W.~Steven Gray}
\address{Old Dominion University, Norfolk, Virginia, 23454, U.S.A.}
         \email{sgray@odu.edu}

\author{Dominique Manchon}
\address{Universit\'e Blaise Pascal,
         C.N.R.S.-UMR 6620, BP 80026,
         63171 Aubi\`ere, France}
         \email{manchon@math.univ-bpclermont.fr}
         \urladdr{http://math.univ-bpclermont.fr/~manchon/}

\thispagestyle{empty}

\date{\today}

\maketitle

\begin{abstract}
The goal of the paper is to give a systematic way to numerically evaluate the generating function of a periodic multiple polylogarithm using a Chen--Fliess series with a rational generating series. The idea is to realize the corresponding Chen--Fliess series as a bilinear dynamical system. A standard form for such a realization is given. The method is also generalized to the case where the multiple polylogarithm has non-periodic components. This allows one, for instance, to numerically validate the Hoffman conjecture.
Finally, a setting in terms of dendriform algebras is provided.
\end{abstract}

\noindent {\footnotesize{{\bf{Keywords}}: multiple polylogarithms, multiple zeta values, Chen--Fliess series, rational formal power series}}

\bigskip

\noindent {\small{{\bf{Math.~Subject Classification}}: 11G55, 11M32, 93B20, 93C10.}}

\tableofcontents


\section{Introduction}
\label{sect:intro}

Given any vector $\mathbf{s}=(s_1,s_2,\ldots,s_l) \in \mathbb{Z}^l$ with $s_1 \geq 2$ and $s_i \geq 1$ for $i \geq 2$, the
associated {\em multiple polylogarithm} (MPL) of {\it{depth}} $l$ and {\it{weight}} $\abs{\mathbf{s}}:=\sum_{i=1}^l s_i$ is taken to be
\begeq \label{eq:Lis}
	\mathrm{Li}_{\mathbf s}(t):=\sum_{k_1>k_2>\cdots >k_l\geq 1} \frac{t^{k_1}}{k_1^{s_1}k_2^{s_2}\cdots k_l^{s_l}},
	\quad \abs{t}\leq 1,
\endeq
whereupon the {\it{multiple zeta value}} (MVZ) of depth $l$ and weight $\abs{\mathbf{s}}$ is the value of \eqref{eq:Lis} at $t=1$, namely,
$$
	\zeta(\mathbf{s}):=\mathrm{Li}_{\mathbf{s}}(1).
$$
Any such vector $\mathbf{s}$ will be referred to as {\em admissible}. The MPL in \eqref{eq:Lis} can be represented in terms of iterated Chen integrals with respect to the $1$-forms $\omega_j^{(1)}:=\frac{dt_j}{1-t_j}$ and $\omega_j^{(0)}:=\frac{dt_j}{t_j}$. Indeed, using the standard notation, $|\mathbf{s}_{(j)}|:=s_1+\cdots+s_j$, $j\in\{1,\ldots,l\}$, one can show that
\begin{equation}
\label{eq:MPLint}
         \mathrm{Li}_{\mathbf s}(t)=\int_{0}^{t} \bigg( \prod_{j=1}^{|\mathbf{s}_{(1)}|-1}
                                \omega_j^{(0)}\bigg) \omega_{|\mathbf{s}_{(1)}|}^{(1)}
                                \cdots \bigg( \prod_{j=|\mathbf{s}_{(l-1)}|+1}^{|\mathbf{s}_{(l)}|-1}
                                \omega_j^{(0)}\bigg) \omega_{|\mathbf{s}_{(l)}|}^{(1)} .
\end{equation}
For instance,
\begin{align*}
	\mathrm{Li}_{(2,1,1)}(t)&=\int_{0}^{t} \frac{dt_1}{t_1}\int_{0}^{t_1} \frac{dt_2}{1-t_2}
								\int_{0}^{t_2} \frac{dt_3}{1-t_3}\int_{0}^{t_3} \frac{dt_4}{1-t_4}
					      = \sum_{k_1>k_2>k_3\geq 1} \frac{t^{k_1}}{k_1^{2}k_2 k_3}.		
\end{align*}
An MPL of depth $l$ is said to be {\em periodic} if it can be written in the form $\mathrm{Li}_{\mathbf{\{s\}}^n}(t)$, where $\{\mathbf{s}\}^n$ denotes the $n$-tuple $(\mathbf{s},\mathbf{s},\ldots,\mathbf{s})\in \mathbb{Z}^{nl}$, $n\geq 0$ with $\mathrm{Li}_{\{\mathbf{s}\}^0}(t):=1$.\footnote{Following other authors, $\{\mathbf{s}\}^n=\{(s_1,s_2,\ldots,s_l)\}^n$ will be written more concisely as $\{s_1,s_2,\ldots,s_l\}^n$.} In this case, the sequence $(\mathrm{Li}_{\{\mathbf{s}\}^n}(t))_{n \in \mathbb{N}_0}$ has the generating function
\begeq
\label{eq:Ls-generating-function}
	\mathcal{L}_{\bf s}(t,\theta):=\sum_{n=0}^\infty \mathrm{Li}_{\{\mathbf{s}\}^n}(t) \left(\theta^{\abs{\mathbf{s}}}\right)^n.
\endeq
In general, the integral representation \eqref{eq:MPLint} implies that $\mathcal{L}_{\mathbf{s}}$ will satisfy a linear ordinary
differential equation in $t$ whose solution can be written in terms of a hypergeometric function \cite{Aoki-etal_07,Borwein-etal_98b,Bowman-Bradley_01,Zakezewski-Zoladek_10,Zakezewski-Zoladek_11,Zakezewski-Zoladek_12,Zoladek_03}.
For example, when $l=1$ and $\mathbf{s}=(s)$, it follows  that
\begeq
\label{eq:L2-ode}
	\left(\left((1-t)\frac{d}{dt}\right)\left(t\frac{d}{dt}\right)^{s-1}-\theta^s\right)\mathcal{L}_{\mathbf{s}}(t,\theta)=0,
\endeq
and its solution is the Euler--Gauss hypergeometric function
\begdi
\label{eq:GF-polylog-L2}
	\mathcal{L}_{(s)}(t,\theta)=_s\! F_{s-1}
\left(\left.
\begin{array}{c}
	-\omega\theta,-\omega^3\theta,\ldots,-\omega^{2s-1}\theta \\
	1,1,\ldots,1
\end{array}
\right|t\right),
\enddi
where $\omega=\expup^{\pi i/s}$, a primitive $s$-th root of $-1$ \cite{Borwein-etal_98b}. By expanding this solution into a hypergeometric series and equating like powers of $\theta$ with those in \rref{eq:Ls-generating-function}, it is possible to show, for example, when $s=2$ that
\begeq
\label{eq:zeta2n-equals-ratio}
	\zeta(\{2\}^n)=\frac{\pi^{2n}}{(2n+1)!},\quad n\geq 1.
\endeq
In a similar manner it can be shown that
\begdi
\label{eq:zeta31n-equals-ratio}
	\zeta(\{3,1\}^n)=\frac{2\pi^{4n}}{(4n+2)!},\quad n\geq 1.
\enddi
This method has yielded a plethora of such MZV identities \cite{Borwein-etal_98a,Borwein-etal_98b,Borwein-Zudilin_xx,Zudilin_03}.
The most general case is treated in \cite{Zoladek_03}, where it is shown that $\mathcal{L}_{\mathbf{s}}$ satisfies the linear differential equation of Fuchs type
\begeq
\label{eq:Ls-odu-Fuch-type}
	(P_{\mathbf{s}}-\theta^{\abs{\mathbf s}})\mathcal{L}_{\mathbf{s}}(t,\theta)=0,
\endeq
where for $\mathbf{s}=(s_1,s_2,\ldots,s_l) \in \mathbb{Z}^l$ 
\begdi
	P_{\mathbf{s}}:=P_{s_l}P_{s_{l-1}}\cdots P_{s_1}
\enddi
and
\begdi
\label{eq:P_si}
	P_{s_i}:=\left((1-t)\frac{d}{dt}\right)\left(t\frac{d}{dt}\right)^{s_i-1}.
\enddi
(The conventions in \cite{Zoladek_03} are to use $-\theta$ in place of $\theta$ and $t$ in place of $1-t$.) In \cite{Zoladek_03} and related work \cite{Zakezewski-Zoladek_10,Zakezewski-Zoladek_11,Zakezewski-Zoladek_12}, the authors develop WKB type asymptotic expansions of these hypergeometric solutions.

\medskip

The ultimate goal of the present paper is to provide a numerical scheme for estimating $\mathcal{L}_{\bf s}(t,\theta)$ by in essence mapping the $\abs{{\bf s}}$-order linear differential equation \rref{eq:Ls-odu-Fuch-type} to a {\em system} of $\abs{{\bf s}}$ first-order bilinear
differential equations which can be solved by standard tools found in software packages like MatLab. Specifically, it will be shown how to construct a dynamical system of the form
\begin{subequations}
\label{eq:general-SISO-bilinear-system}
\begin{align}
	\dot{z}&=N_0z\,u_0+N_1z\,u_1,\;\;z(0)=z_0 \label{eq:bilinear-SISO-state-equation}\\
		y&=Cz,
\end{align}
\end{subequations}
which when simulated over the interval $[0,1]$ has the property that $y(t)=\mathcal{L}_{\mathbf{s}}(t,\theta)$ for any value of $\theta$ and $t\in[0,1]$. In this case, the matrices $N_0$ and $N_1$ will depend on $\theta$, and the initial condition $z_0$ and the input functions $u_0$, $u_1$ must be suitably chosen. Such a technique could be useful for either disproving certain conjectures involving MZVs or providing additional evidence for the truthfulness of other conjectures. For example, one could validate with a certain level of (numerical) confidence a conjecture of the form
\begdi
	\zeta(\{\mathbf{s}_a\}^n)=b^n\zeta(\{\mathbf{s}_b\}^n),\quad n \in \mathbb{N}, b \in \mathbb{Z},
\enddi
where $\mathbf{s}_a\in\mathbb{Z}^{l_a}$, $\mathbf{s}_b\in\mathbb{Z}^{l_b}$ with $\abs{\mathbf{s}_a}=\abs{\mathbf{s}_b}$.
Take as a specific example the known identity
\begeq
\label{eq:zeta4n-equals-4zeta31n}
	\zeta(\{4\}^n)=4^n\zeta(\{3,1\}^n)%
\endeq
for all $n \ge 1$, so that $\mathbf{s}_a=(4)$, $\mathbf{s}_b=(3,1)$ and $b=4$ \cite{Borwein-etal_98b}. Note that for $n=1$ the identity follows immediately from double shuffle relations for MZVs \cite{IKZ_06}. On the level of generating functions it is evident that
\begin{align*}
	\mathcal{L}_{(4)}(1,\theta)&=\sum_{n=0}^\infty \mathrm{Li}_{\{4\}^n}(1) \left(\theta^{4}\right)^n
					=\sum_{n=0}^\infty \zeta(\{4\}^n)\,\theta^{4n} \\
	\mathcal{L}_{(3,1)}(1,\sqrt{2}\theta)&=\sum_{n=0}^\infty \mathrm{Li}_{\{3,1\}^n}(1) \left((\sqrt{2}\theta)^{4}\right)^n
					=\sum_{n=0}^\infty 4^n\zeta(\{3,1\}^n)\,\theta^{4n}.
\end{align*}
Therefore, identity \rref{eq:zeta4n-equals-4zeta31n} implies that
\begeq
\label{eq:L4-L31}
	\mathcal{L}_{(4)}(1,\theta)-\mathcal{L}_{(3,1)}(1,\sqrt{2}\theta)=0,\quad \forall\theta\in \re,
\endeq
a claim that can be tested empirically if these generating functions can be accurately evaluated. The method can also be generalized to address the conjecture of Hoffman that
\begeq
\label{eq:Hoffman-conjecture}
	\zeta(\{2\}^{n},2,2,2)+2\zeta(\{2\}^n,3,3)=\zeta(2,1,\{2\}^n,3),
\endeq
for all integers $n>0$, which has only been proved for $n\leq 8$ \cite{Borwein-Zudilin_xx}. The idea here is to admit {\em non-periodic} components in the generating function calculation. For example, $(\{2\}^n,3,3)$ can be viewed as having the periodic component $\{2\}^n$ and the non-periodic component $(3,3)$. In the general case, say when $\mbf{s}_n:=(\mathbf{s}_a,\{\mathbf{s}_b\}^n,\mathbf{s}_c)$, $n\geq 0$, the generating function is defined analogously as
\begdi
	\mathcal{L}_{(\mathbf{s}_a,\{\mathbf{s}_b\},\mathbf{s}_c)}(t,\theta):=
		\sum_{n=0}^\infty {\rm Li}_{{\mathbf s}_n}(t)\left(\theta^{\abs{\mathbf{s}_b}}\right)^n.
\enddi
Therefore, relation \rref{eq:Hoffman-conjecture}, if true, would imply that
\begeq
\label{eq:GF-identity-Hoffman-conjecture}
	\mathcal{L}_{(\{2\},2,2,2)}(1,\theta)+2\mathcal{L}_{(\{2\},3,3)}(1,\theta)-\mathcal{L}_{(2,1,\{2\},3)}(1,\theta)
	=0,\quad \forall\theta\in\re.
\endeq
The basic approach to estimating $\mathcal{L}_{\bf s}(t,\theta)$ is to map a periodic multiple polylogarithm to a rational series and then to employ concepts from control theory to produce bilinear state space realization \rref{eq:general-SISO-bilinear-system} of the corresponding rational Chen--Fliess series \cite{Berstel-Reutenauer_10,Elliott_09,Fliess_81}. The periodic nature of the MPL always ensures that these realizations have a certain built-in recursion/feedback structure. The technique will first be described in general, and then it will be demonstrated by empirically verifying the identities \rref{eq:zeta2n-equals-ratio}, \rref{eq:zeta4n-equals-4zeta31n}, and \rref{eq:Hoffman-conjecture}.

\medskip

The paper is organized as follows. In the next section, a brief summary of rational Chen--Fliess series is given to establish the notation and the basic concepts to be employed. Then the general method for evaluating a generating function of a periodic multiple polylogarithm is given in the subsequent section, which also contains in Subsection \ref{ssect:dend} a short digression regarding another way of looking at periodic MPLs in terms of the shuffle algebra. This is followed by several examples in Section \ref{sec:examples}. In particular, the last example shows that the Hoffman conjecture \rref{eq:Hoffman-conjecture} has a high likelihood of being true. The final section gives the paper's conclusions.


\section{Preliminaries}


\subsection{Chen--Fliess series}
\label{ssect:ChFseries}

A finite nonempty set of noncommuting symbols $X=\{ x_0,x_1,$ $\ldots,x_m\}$ is called an {\em alphabet}. Each element of $X$ is called a {\em letter}, and any finite sequence of letters from $X$, $\eta=x_{i_1}\cdots x_{i_k}$, is called a {\em word} over $X$. The {\em length} of word $\eta$, denoted $\abs{\eta}$, is the number of letters in $\eta$. The set of all words with fixed length $k$ is denoted by $X^k$. The set of all words including the empty word, $\emptyset$, is designated by $X^\ast$. It forms a monoid under catenation.
The set $\eta X^\ast\xi\subseteq X^\ast$ is the set of all words with prefix $\eta$ and suffix $\xi$.
Any mapping $c: X^\ast\rightarrow \re^\ell$ is called a {\em formal power series}. The value of $c$ at $\eta \in X^\ast$ is written as $(c,\eta)\in \re^\ell$ and called the {\em coefficient} of the word $\eta$ in the series $c$. Typically, $c$ is  represented as the formal sum $c=\sum_{\eta\in X^\ast}(c,\eta)\eta.$ If the {\em constant term} $(c,\emptyset)=0$ then $c$ is said to be {\em proper}. The collection of all formal power series over the alphabet $X$ is denoted by $\allseriesell$. The subset of polynomials is written as $\allpolyell$. Each set forms an associative $\re$-algebra under the catenation product.

\begin{defn}\label{def:left-shift operator}
Given $\xi \in X^{\ast}$, the corresponding \bfem{left-shift operator} $\xi^{-1}:X^{\ast} \rightarrow X^{\ast}$ is defined:
\begdi
 \eta \mapsto \xi^{-1}(\eta):=\left\{ \begar{ccl}
                     \eta^{\prime} & : &\eta=\xi \eta^{\prime} \\[0.03in]
                      0 & : & \hbox{otherwise}.
                          \endar\right.
\enddi
It is extended linearly to $\allseriesell$.
\end{defn}
 One can formally associate with any series $c\in\allseriesell$ a causal $m$-input, $\ell$-output operator, $F_c$, in the following manner. Let $\mathfrak{p} \ge 1$ and $t_0 < t_1$ be given. For a Lebesgue measurable function $u: [t_0,t_1] \rightarrow\re^m$, define $\norm{u}_{\mathfrak{p}} := \max\{\norm{u_i}_{\mathfrak{p}}: \ 1\le i \le m\}$, where $\norm{u_i}_{\mathfrak{p}}$ is the usual $L_{\mathfrak{p}}$-norm for a measurable real-valued function, $u_i$, defined on the interval  $[t_0,t_1]$. Let $L^m_{\mathfrak{p}}[t_0,t_1]$ denote the set of all measurable functions defined on $[t_0,t_1]$ having a finite $\norm{\cdot}_{\mathfrak{p}}$ norm and $B_{\mathfrak{p}}^m(R)[t_0,t_1]:=\{u\in L_{\mathfrak{p}}^m[t_0,t_1]:\norm{u}_{\mathfrak{p}}\leq R\}$. Assume $C[t_0,t_1]$ is the subset of continuous functions in $L_{1}^m[t_0,t_1]$. Define inductively for each word $\eta \in X^{\ast}$ the map $E_\eta: L_1^m[t_0, t_1] \rightarrow C[t_0, t_1]$ by setting $E_\emptyset[u]=1$ and letting
\[
	E_{x_i\bar{\eta}}[u](t,t_0) := \int_{t_0}^tu_{i}(\tau)E_{\bar{\eta}}[u](\tau,t_0)\,d\tau,
\]
where $x_i\in X$, $\bar{\eta}\in X^{\ast}$, and $u_0=1$. The input-output operator corresponding to the series $c \in \allseriesell$ is the {\em Fliess operator} or {\em Chen--Fliess series}
\begeq
\label{eq:Fliess-operator-defined}
	F_c[u](t) = \sum_{\eta\in X^{\ast}} (c,\eta)\,E_\eta[u](t,t_0).
\endeq
If there exist real numbers $K_c,M_c>0$ and $r\in[0,1)$ such that the coefficients of the generating series  $c=\sum_{\eta\in X^\ast}(c,\eta)\eta \in \allseriesell$ satisfying
the growth bound
\begeq
	\abs{(c,\eta)}\le K_c M_c^{|\eta|}(\abs{\eta}!)^r,\quad \forall\eta\in X^{\ast}, \label{eq:global-convergence-growth-bound}
\endeq
then the series \rref{eq:Fliess-operator-defined} defines an operator from the extended space $\Lpme (t_0)$ into $C[t_0, \infty)$,  where
\begdi
	\Lpme(t_0) := \{u:[t_0,\infty)\rightarrow \re^m:u_{[t_0,t_1]}\in \Lpm[t_0,t_1],\; \forall t_1 \in (t_0,\infty)\},
\enddi
and $u_{[t_0,t_1]}$ denotes the restriction of $u$ to the intervall $[t_0,t_1]$. (Here, $\abs{z}:=\max_i \abs{z_i}$ when $z\in\re^\ell$.) See \cite{Winter-Gray-Duffaut_Espinosa_CISS15} for details. In this case, the operator is said to be {\em globally convergent}, and the set of all series satisfying \rref{eq:global-convergence-growth-bound} is designated by $\allseriesellGC$. In the following sections, it suffices to set $\ell=m=1$, which corresponds to the single-input, single-output (SISO) case.


\subsection{Bilinear realizations of rational Chen--Fliess series}
\label{ssect:BilinReal}

A series $c\in\allseries$ is called {\em invertible} if there exists a series $c^{-1}\in\allseries$ such that $cc^{-1}=c^{-1}c=1$.\footnote{The polynomial $1\emptyset$ is abbreviated throughout as $1$.} In the event that $c$ is not proper, i.e., the coefficient $(c,\emptyset)$ is nonzero, it is always possible to write
\begdi
	c=(c,\emptyset)(1-c^{\prime}),
\enddi
where $c^{\prime}\in\allseries$ is proper.
It then follows that
\begdi
	c^{-1}=\frac{1}{(c,\emptyset)}(1-c^{\prime})^{-1}
		 =\frac{1}{(c,\emptyset)} (c^\prime)^{\ast},
\enddi
where
\begdi
	(c^{\prime})^{\ast}:=\sum_{i=0}^{\infty} (c^{\prime})^i.
\enddi
In fact, $c \in\allseries$ is invertible {\em if and only if} $c$ is not proper. Now let $S$ be a subalgebra of the $\re$-algebra $\allseries$ with the catenation product. $S$ is said to be {\em rationally closed} when every invertible $c\in S$ has $c^{-1}\in S$ (or equivalently, every proper $c^{\prime}\in S$ has $(c^{\prime})^{\ast}\in S$). The {\em rational closure} of any subset $E\subset\allseries$ is the smallest rationally closed subalgebra of $\allseries$ containing $E$.

\begde
A series $c\in\allseries$ is \bfem{rational} if it belongs to the rational closure of $\allpoly$.
\endde

It turns out that an entirely different characterization of a rational series is possible using the following concept.

\begde \label{def:linear-representation}
A \bfem {linear} \bfem{representation} of a series $c\in\allseries$ is any triple $(\mu,\gamma,\lambda)$, where
\begdi
	\mu:X^{\ast}\rightarrow \re^{n\times n}
\enddi
is a monoid morphism, and the vectors $\gamma,\lambda^T\in\re^{n\times 1}$ are such that each coefficient
\begdi
	(c,\eta)=\lambda\mu(\eta)\gamma, \quad \forall\eta\in X^{\ast}.
\enddi
The integer $n$ is the dimension of the representation.
\endde

\begde
A series $c\in\allseries$ is called \bfem{recognizable} if it has a linear representation.
\endde

\begth {\rm (Sch\"{u}tzenberger)} \label{th:Schutzenberger}
A formal power series is rational if and only if it is recognizable.
\endth

 Returning to \eqref{eq:Fliess-operator-defined}, Chen--Fliess series $F_c$ is said to be rational when its generating series $c\in\allseries$ is rational. The state space realization \rref{eq:general-SISO-bilinear-system} is said to {\em realize $F_c$ on $\Lpe(t_0)$} when (\ref{eq:bilinear-SISO-state-equation}) has a well defined solution, $z(t)$, on the interval $[t_0,t_0+T]$ for every $T>0$ with input $u\in \Lpe(t_0)$  and output
\begdi
	y(t)=F_c[u](t)=C(z(t)),\quad t\in[t_0,t_0+T].
\enddi
Identify with any linear representation $(\mu,\gamma,\lambda)$ of the series $c\in\allseries$ the bilinear system
$$
	(N_0,N_1,z_0,C):=(\mu(x_0),\mu(x_1),\gamma,\lambda).
$$
The following result is well known \cite{Fliess_81,Gray-Wang_SCL02}.

\begth \label{th:finite-Hankel-rank-augmented}
The statements below are equivalent for a given $c\in\allseries$:
\begin{description}
	\item[i] $(\mu,\gamma,\lambda)$ is a linear representation of $c$.
	\item[ii] The bilinear system $(N_0,N_1,z_0,C)$ realizes $F_c$ on $\Lpe(t_0)$ for any $\mathfrak{p}\geq 1$.
\end{description}
\endth



\section{Evaluating periodic multiple polylogarithms}
\label{sect:PeriodicMPL}

It is first necessary to associate a periodic MPL and its generating function to a rational series. Elements of this idea have appeared in numerous places. The approach taken here is most closely related to the one presented in \cite{Houseaux-etal_03}. The next step is then to find the bilinear realization of the rational Chen--Fliess series in terms of its linear representation (see Theorem~\ref{th:realization-GF-periodic-multiple-polylogs}). The case when non-periodic components are present works similarly
but is slightly more complicated (see Theorem~\ref{th:realization-GF-nonperiodic-multiple-polylogs}). Recall that throughout $m=1$,
so that the underlying alphabet is $X:=\{x_0,x_1\}$.


\subsection{Periodic multiple polylogarithms}
\label{ssect:PeriodicMPL}

Given any admissible vector $\mathbf{s} \in \mathbb{Z}^l$, there is an associated word $\eta_{\mathbf s} \in x_0X^\ast x_1$ of length $|\mathbf s|$
\begdi
	\eta_{\mathbf s}=x_0^{s_1-1}x_1x_0^{s_2-1}x_1\cdots x_0^{s_l-1}x_1.
\enddi
In which case, $c_{\mathbf{s}}:=(\theta^{\abs{\mathbf{s}}} \eta_{\mathbf{s}})^\ast=\sum_{n\geq 0} \left(\theta^{\abs{\mathbf{s}}}\eta_{\mathbf{s}}\right)^n$ is a rational series satisfying the identity
\begeq
\label{eq:etas-is-periodic}
	1+(\theta^{\abs{\mathbf{s}}}\eta_{\mathbf{s}})c_{\mathbf{s}}=c_{\mathbf{s}}.
\endeq
The idea is to now relate the generating function of the sequence $(\mathrm{Li}_{\{\mathbf{s}\}^n}(t))_{n>0}$ to the Chen--Fliess series with generating series $c_{\mathbf{s}}$. Recall that for any word $x_i\xi^\prime\in X^\ast$ the iterated integral is defined inductively by
\begdi
	E_{x_i\xi^\prime}[u](t)=\int_0^t u_{i}(\tau) E_{\xi^\prime}[u](\tau)\,d\tau,
\enddi
where $x_i\in X$, $\xi^\prime\in X^\ast$. Assume here that the letters $x_0$ and $x_1$ correspond to the inputs $u_0(t):=1/t$ and $u_1(t):=1/(1-t)$, respectively,
and $E_{\emptyset}:=1$. For the formal power series $c_{\mathbf s}\in\allseries$, the corresponding Chen--Fliess series is then taken to be
\begdi
	F_{c_{\mathbf s}}[u]=\sum_{\xi\in X^\ast} (c_{\mathbf s},\xi)E_\xi[u].
\enddi
Comparing this to the classical definition \rref{eq:Fliess-operator-defined}, the factor $1/t$ can be extracted from $u_0$ and $u_1$ so that each integral can be viewed instead as integration with respect to the Haar measure. That is,
\begdi
	E_{x_i\xi^\prime}[u](t)=\int_0^t \bar{u}_{i}(\tau) E_{\xi^\prime}[u](\tau)\frac{d\tau}{\tau},
\enddi
where $\bar{u}_0(t):=1$ and $\bar{u}_1(t)=tu_1(t)$. The following lemma now applies.

\begin{lem} \label{le:Ls-equals-Fcs}
For any admissible vector $\mathbf{s}\in\mathbb{Z}^l$,
\begdi
	\mathcal{L}_{\mathbf{s}}(t,\theta)=F_{c_{\mathbf{s}}}[\mathrm{Li}_0](t),\;\;t\in[0,1],\;\; \theta\in\re,
\enddi
where $\mathrm{Li}_0(t):=t/1-t$.
\end{lem}

\begpr
First observe that since $c_{\mathbf{s}}=\sum_{n\geq 0} \left(\theta^{\abs{\mathbf{s}}}\eta_{\mathbf{s}}\right)^n$, it follows directly that
\begdi
	F_{c_{\mathbf{s}}}[u](t)	=\sum_{n=0}^\infty F_{\left(\theta^{\abs{\mathbf{s}}}\eta_{\mathbf{s}}\right)^n}[u](t)
						=\sum_{n=0}^\infty E_{\eta_{\mathbf{s}}^n}[u](t)\left(\theta^{\abs{\mathbf{s}}}\right)^n.
\enddi
Comparing this against the definition
\begdi
	\mathcal{L}_{\bf s}(t,\theta)=\sum_{n=0}^\infty \mathrm{Li}_{\{\mathbf{s}\}^n}(t) \left(\theta^{\abs{\mathbf{s}}}\right)^n,
\enddi
it is evident that one only needs to verify the identity
\begeq \label{eq:Fetas-equals-Lis}
E_{\eta_{\mathbf{s}}^n}[\mathrm{Li}_0](t)=\mathrm{Li}_{\{\mathbf{s}\}^n}(t),\;\; n\geq 0.
\endeq
But this is clear from \eqref{eq:MPLint}, i.e., for any admissible vector $\mathbf{s} \in \mathbb{Z}^l$
\begdi
	\mathrm{Li}_{\mathbf{s}}(t)=\int_{0}^t u_i(\tau)\mathrm{Li}_{\mathbf{s}^\prime}(\tau)\,d\tau,
\enddi
where $\eta_{\mathbf s}=x_i\eta_{{\mathbf s}^\prime}$,
\begdi
	u_i(t)=\left\{
		\begin{array}{ccl}
			\frac{1}{t}			&:& i=0 \\[0.05in]
			\frac{t}{1-t}\frac{1}{t}	&:& i=1,
		\end{array}
\right.
\enddi
and $\mathrm{Li}_\emptyset(t)=1$ \cite{Zudilin_03}. Therefore, it follows directly that $\mathrm{Li}_{\mathbf{s}}(t)=E_{\eta_{\mathbf s}}[\mathrm{Li}_0](t)$, from which \rref{eq:Fetas-equals-Lis} also follows.
\endpr

The key idea now is to apply Theorem~\ref{th:finite-Hankel-rank-augmented} and the rational nature of the series $c_{\mathbf{s}}$ in order to build a bilinear realization of the mapping $u\mapsto y=F_{c_{\mathbf{s}}}[u]$ so that $\mathcal{L}_{\mathbf{s}}(t,\theta)$ can be evaluated by numerical simulation of a dynamical system. In principle, one could attempt to ensure
 that any such realization is minimal in dimension or even canonical in some sense \cite{Brockett_72,DAlessandro-etal_74,Dorissen_89,Sussman_76}, but in the present context these properties are not really essential.

\begth \label{th:realization-GF-periodic-multiple-polylogs}
For any admissible $\mathbf{s}\in\mathbb{Z}^l$, $\mathcal{L}_{\mathbf{s}}(t,\theta)=F_{c_{\mathbf{s}}}[\mathrm{Li}_0](t)$ has the bilinear realization
$$
	(N_0,N_1,z_0,C):=\big(\mu(x_0),\mu(x_1),\gamma,\lambda\big),
$$
where
\begin{subequations}
\label{eq:bilinear-Fetas}
\begin{align}
	N_0&=\diag\big(N_0(s_1),N_0(s_2),\ldots,N_0(s_l)\big) \\[0.3cm]
N_1&=I^+_{\abs{\bf s}}-N_0+\theta^{\abs{s}}e_{\abs{\bf s}}e_1^T
\end{align}
\end{subequations}
with $N_0(s_i)\in\re^{s_i\times s_i}$ and $I^+_{\abs{\bf s}}\in \re^{\abs{\bf s}\times \abs{\bf s}}$ being matrices of zeros except for a super diagonal of ones, $e_{i}$ is an elementary vector with a one in the $i$-th position, and $z_0=C^T=e_1\in \re^{\abs{\mathbf{s}}\times 1}$.
\endth

\begpr
First recall Definition \ref{def:left-shift operator} describing
the left-shift operator on $X^*$, i.e., for any $x_i \in X$, $x_i^{-1}(\cdot)$ is defined by $x_i^{-1}(x_i\eta)=\eta$ with $\eta\in X^\ast$ and zero otherwise. In which case, $(x_i\xi)^{-1}(\cdot)=\xi^{-1}x_i^{-1}(\cdot)$ for any $\xi\in X^\ast$. Now assign the first state of the realization to be
\begdi
	z_1(t)=F_{c_{\mathbf{s}}}[u](t)=1+F_{(\theta^{\abs{\mathbf{s}}}\eta_{\mathbf{s}})c_{\mathbf{s}}}[u](t).
\enddi
 In light of the integral representation \eqref{eq:MPLint} of MPLs, differentiate $z_1$ exactly $s_1$ times so that the input $u_1(t):=\bar{u}_1(t)/t$ appears. Assign a new state at each step along the way. Specifically,
\begin{align*}
	\dot{z}_1(t)&=\frac{1}{t}F_{\theta^{\abs{\mathbf{s}}}x_0^{-1}(\eta_{\mathbf{s}}) c_{\mathbf{s}}}[u](t)=:z_2(t)\frac{1}{t} \\
			 &\hspace*{0.08in}\vdots \\
	\dot{z}_{s_1-1}(t)&=\frac{1}{t}F_{\theta^{\abs{\mathbf{s}}}(x_0^{s_1-1})^{-1}(\eta_{\mathbf{s}}) c_{\mathbf{s}}}[u](t)
					=:z_{s_1}(t)\frac{1}{t} \\
	\dot{z}_{s_1}(t)&=\bar{u}_1(t)\frac{1}{t}F_{\theta^{\abs{\mathbf{s}}}(x_0^{s_1-1}x_1)^{-1}(\eta_{\mathbf{s}}) c_{\mathbf{s}}}[u](t)
					=:z_{s_1+1}(t)\bar{u}_1(t)\frac{1}{t}.
\end{align*}
This produces the first $s_1$ rows of the matrices  in \rref{eq:bilinear-Fetas}
since when $l>1$
\begin{align*}
\left[
\begin{array}{c}
\dot{z}_1(t) \\
\vdots \\
\dot{z}_{s_1-1}(t) \\
\dot{z}_{s_1}(t)
\end{array}
\right]&=I^+_{s_1\times (s_1+1)}
\left[
\begin{array}{c}
z_1(t) \\
\vdots \\
z_{s_1}(t) \\
z_{s_1+1}(t)\bar{u}_1(t)
\end{array}
\right]\frac{1}{t} \\
&=
\left[
\begin{array}{c|c}
N_0(s_1) & 0
\end{array}
\right]
\left[
\begin{array}{c}
z_1(t) \\
\vdots \\
z_{s_1}(t) \\ \hline
z_{s_1+1}(t)
\end{array}
\right]
\frac{1}{t}+
\left[
\begin{array}{c|c}
{\bf 0}_{s_1} & e_{s_1}
\end{array}
\right]
\left[
\begin{array}{c}
z_1(t) \\
\vdots \\
z_{s_1}(t) \\ \hline
z_{s_1+1}(t)
\end{array}
\right]
\bar{u}_1(t)\frac{1}{t}.
\end{align*}%
Both $\left[\begin{array}{c|c}
N_0(s_1) & 0
\end{array}
\right]$
and
$\left[
\begin{array}{c|c}
{\bf 0}_{s_1} & e_{s_1}
\end{array}
\right]$
denote matrices in $\re^{s_1\times (s_1+1)}$.  The pattern is exactly repeated until the final state, then the periodicity of $c_{\mathbf{s}}$ comes into play. Namely,
\begdi
	\dot{z}_{\abs{\mathbf s}}(t)= \theta^{\abs{\mathbf{s}}}\bar{u}_1(t)\frac{1}{t}F_{(\eta_{\mathbf{s}})^{-1}(\eta_{\mathbf{s}}) c_{\mathbf{s}}}[u](t)
	=:\theta^{\abs{\mathbf{s}}}z_{1}(t)\bar{u}_1(t)\frac{1}{t},
\enddi
which gives the final rows of $N_0$ and $N_1$ in \rref{eq:bilinear-Fetas}.
\endpr

It is worth pointing out that the validity of \rref{eq:Ls-odu-Fuch-type} is obvious in the present setting. Namely, \rref{eq:Ls-odu-Fuch-type} follows from the fact that \rref{eq:etas-is-periodic} implies $\eta_{\mathbf s}^{-1}(c_{\mathbf s})-\theta^{\abs{\mathbf{s}}}c_{\mathbf s}=0$, and thus, Lemma~\ref{le:Ls-equals-Fcs} gives
\begdi
	(P_{\mathbf s}-\theta^{\abs{\mathbf s}})\mathcal{L}_{\mathbf{s}}(t,\theta)
	= (P_{\mathbf s}-\theta^{\abs{\mathbf s}})F_{c_{\mathbf s}}[\mathrm{Li}_0](t)
	= F_{\eta_{\mathbf s}^{-1}(c_{\mathbf s})-\theta^{\abs{\mathbf{s}}}c_{\mathbf s}}[\mathrm{Li}_0](t)
	= F_{0\cdot c_{\mathbf s}}[\mathrm{Li}_0](t)=0.
\enddi


\subsection{Periodic multiple polylogarithms with non-periodic components}
\label{ssect:NonStrictPeriodicMPL}

The non-periodic case requires a generalization of the basic set-up.  The following lemma links this class of generating functions
to the corresponding set of rational Fliess operators.

\begle \label{le:GF-nonperiodic-multiple-polylogs}
For any admissible $\mathbf{s}:=(\mathbf{s}_a,\{\mathbf{s}_b\},\mathbf{s}_c)$
\begdi
	\mathcal{L}_{\mathbf{s}}(t,\theta)=F_{c_{\mathbf{s}}}[{\rm Li}_0](t),\;\; t\in[0,1],\;\;\theta\in\re,
\enddi
where $c_{\mathbf{s}}:=\eta_{\mathbf{s}_a}\left(\theta^{\abs{\mathbf{s}_b}}\eta_{\mathbf{s}_b}\right)^\ast \eta_{\mathbf{s}_c}$.
\endle

\begpr
Similar to the periodic case, $c_{\mathbf{s}}=\sum_{n\geq 0} \eta_{{\bf s}_a}\left(\theta^{\abs{\mathbf{s}_b}}\eta_{\mathbf{s}_b}\right)^n\eta_{{\bf s}_c}$, and therefore,
\begin{align*}
	F_{c_{\mathbf{s}}}[u](t)&=\sum_{n=0}^\infty F_{\eta_{\mathbf{s}_a}\left(\theta^{\abs{\mathbf{s}_b}}
		\eta_{\mathbf{s}_b}\right)^n\eta_{\mathbf{s}_c}}[u](t)
	=\sum_{n=0}^\infty E_{\eta_{\mathbf{s}_a}\eta_{\mathbf{s}_b}^n\eta_{\mathbf{s}_c}}[u](t)\left(\theta^{\abs{\mathbf{s}_b}}\right)^n.
\end{align*}
The same argument used for proving \rref{eq:Fetas-equals-Lis} now shows that $E_{\eta_{\mathbf{s}_a}\eta_{\mathbf{s}_b}^n\eta_{\mathbf{s}_c}}[\mathrm{Li}_0](t)= {\rm Li}_{{\mathbf s}_n}(t)$, $n\geq 0$. In which case, $F_{c_{\mathbf{s}}}[\mathrm{Li}_0](t)=\mathcal{L}_{\mathbf{s}}(t,\theta)$ as claimed.
\endpr

The required generalization of Theorem~\ref{th:realization-GF-periodic-multiple-polylogs} is a bit more complicated. A simple example is given first to motivate the general approach.

\begex \label{ex:realization-L212n3}
Consider the periodic MPL with non-periodic components specified by $\mbf{s}=(2,1,\{2\},3)$ as
appearing in \rref{eq:GF-identity-Hoffman-conjecture}. In this case, $c_{\mathbf s}=\sum_{n\geq 0}x_0x_1^2(\theta^2x_0x_1)^nx_0^2x_1=x_0x_1^2\bar{c}$, where $\bar{c}=x_0^2x_1+\theta^2x_0x_1\bar{c}$. Assign the first state of the realization to be
\begdi
	z_1(t)=F_{c_{\mathbf{s}}}[u](t)=F_{x_0x_1^2\bar{c}}[u](t).
\enddi
The strategy here is to differentiate $z_1$ exactly $\abs{\eta_{\mathbf{s}_a}}=\abs{x_0x_1^2}=3$ times, assigning new states along the way, in order to remove the prefix $x_0x_1^2$ and isolate $\bar{c}$. At which point, the identity $\bar{c}=x_0^2x_1+\theta^2x_0x_1\bar{c}$ is used and the process is continued. This will yield a certain block diagonal structure for $N_0$ and an upper triangular form for $N_1$. As will be shown shortly, this structure is completely general but possibly redundant. Specifically,
\begin{align*}
	\dot{z}_1(t)&= \frac{1}{t}F_{x_1^2\bar{c}}[u](t)
				=:z_2(t)\frac{1}{t} \\
	\dot{z}_2(t)&=\frac{1}{t}\bar{u}_1(t)F_{x_1\bar{c}}[u](t)
				=:z_3(t)\bar{u}_1(t)\frac{1}{t} \\
	\dot{z}_3(t)&=\frac{1}{t}\bar{u}_1(t)F_{\bar{c}}[u](t)
				=\frac{1}{t}\bar{u}_1(t)F_{x_0^2x_1+\theta^2x_0x_1\bar{c}}[u](t)
				=:z_4(t)\bar{u}_1(t)\frac{1}{t} \\
	\dot{z}_4(t)&=\frac{1}{t}F_{x_0x_1+\theta^2x_1\bar{c}}[u](t)
				=:z_5(t)\frac{1}{t} \\
	\dot{z}_5(t)&=\frac{1}{t}F_{x_1}[u](t)+\frac{\theta^2}{t}\bar{u}_1(t)F_{\bar{c}}[u](t)
				=:z_6(t)\frac{1}{t}+\theta^2z_4(t)\bar{u}_1(t)\frac{1}{t} \\
	\dot{z}_6(t)&=\bar{u}_1(t)\frac{1}{t}.
\end{align*}
The corresponding realization at this point has the form
\begin{align*}
\dot{z}=&\tilde{N}_0z \bar{u}_0+\tilde{N}_1z\bar{u}_1+B_1\bar{u}_1,\;\;z(0)=\tilde{z}_0 \\
y=&\tilde{C}z,
\end{align*}
which does not have the form of a bilinear realization as defined in \rref{eq:general-SISO-bilinear-system} since the state equation for
$z_6$ does not depend on $z$, and thus, the term $B_1\bar{u}_1$ with $B_1=e_6$ appears.
Nevertheless, a permutation of the canonical embedding of Brockett (see \cite[Theorem 1]{Brockett_72}), namely,
\begeq \label{eq:brockett-embedding}
N_0=
\left[
\begin{array}{cc}
\tilde{N}_0 & 0 \\
0 & 0
\end{array}
\right],\;\;
N_1=
\left[
\begin{array}{cc}
\tilde{N}_1 & B_1 \\
0 & 0
\end{array}
\right],\;\;
z_0=
\left[
\begin{array}{c}
\tilde{z}_0 \\
1
\end{array}
\right],\;\;
C^T=
\left[
\begin{array}{c}
\tilde{C}^T \\
0
\end{array}
\right],
\endeq
renders an input-output equivalent bilinear realization of the desired form. In this case,
\begdi
N_0=
\left[
\begin{array}{ccc|ccc|c}
0 & 1 & 0 & 0 & 0 & 0 & 0 \\
0 & 0 & 0 & 0 & 0 & 0 & 0 \\
0 & 0 & 0 & 0 & 0 & 0 & 0\\ \hline
0 & 0 & 0 & 0 & 1 & 0 & 0\\
0 & 0 & 0 & 0 & 0 & 1 & 0\\
0 & 0 & 0 & 0 & 0 & 0 & 0\\ \hline
0 & 0 & 0 & 0 & 0 & 0 & 0\\
\end{array}
\right],\;\;
N_1=
\left[
\begin{array}{ccc|ccc|c}
0 & 0 & 0 & 0 & 0 & 0 & 0 \\
0 & 0 & 1 & 0 & 0 & 0 & 0 \\
0 & 0 & 0 & 1 & 0 & 0 & 0 \\ \hline
0 & 0 & 0 & 0 & 0 & 0 & 0 \\
0 & 0 & 0 & \theta^2 & 0 & 0 & 0 \\
0 & 0 & 0 & 0 & 0 & 0 & 1 \\ \hline
0 & 0 & 0 & 0 & 0 & 0 & 0 \\
\end{array}
\right],\;\;
z(0)=\left[\begin{array}{c}
1 \\
0 \\
0 \\ \hline
0 \\
0 \\
0 \\ \hline
1 \\
\end{array}\right],\;\;
C^T=
\left[\begin{array}{c}
1 \\
0 \\
0 \\ \hline
0 \\
0 \\
0 \\ \hline
0
\end{array}\right].
\enddi
\endex
\begth \label{th:realization-GF-nonperiodic-multiple-polylogs}
Consider any admissible $\mathbf{s}:=(\mathbf{s}_a,\{\mathbf{s}_b\},\mathbf{s}_c)$ with $\eta_{\mathbf{s}_a}:=x_{i_1}\cdots x_{i_k}$, $k={j_{\abs{\mathbf{s}_a}}}$, and $\abs{\mathbf{s}_c}>0$. Then $\mathcal{L}_{\mathbf{s}}(t,\theta)=F_{c_{\mathbf{s}}}[\mathrm{Li}_0](t)$
has the bilinear realization $(N_0,N_1,z_0,C)$, where
\begdi
	N_0=\diag(N_0(\mathbf{s}_a),N_0(\mathbf{s}_b,\mathbf{s}_c),0),\;\;
	N_1=\left[
\begin{array}{cc}
	N_1(\mathbf{s}_a)	& E_{\abs{\mathbf{s}_a} 1}\\
				0	& N_1(\mathbf{s}_b,\mathbf{s}_c)
\end{array}
\right]
\enddi
with $N_i(\mathbf{s}_a)\in\re^{\abs{\mathbf{s}_a}\times \abs{\mathbf{s}_a}}$ being a matrix of zeros and ones depending only on $\mathbf{s}_a$, $E_{\abs{\mathbf{s}_a} 1}$ is the elementary matrix with a one in position $(\abs{\mathbf{s}_a},1)$, and $N_i(\mathbf{s}_b,\mathbf{s}_c)\in\re^{s_{bc}\times s_{bc}}$ is a matrix of zeros, ones, and the entry $\theta^{\abs{\mathbf{s}_b}}$. (Its dimension $s_{bc}$ and exact structure depend only on $\mathbf{s}_b$ and $\mathbf{s}_c$.) Finally,
$z_0=e_1+e_{\abs{\mathbf{s}_a}+s_{bc}}\in \re^{(\abs{\mathbf{s}_a}+s_{bc})\times 1}$
and $C=e_1\in \re^{1\times (\abs{\mathbf{s}_a}+s_{bc})}$.
\endth

\begpr
Following Example~\ref{ex:realization-L212n3}, assign the first state of the realization to be
\begdi
	z_1(t)=F_{c_{\mathbf{s}}}[u](t)
		 =F_{\eta_{\mathbf{s}_a}\bar{c}}[u](t),
\enddi
where $\bar{c}:=\eta_{{\mathbf s}_c}+\theta^{\abs{\eta_{\mathbf{s}_b}}}\eta_{\mathbf{s}_b}\bar{c}$,
and differentiate $z_1$ until the series $\bar{c}$ appears in isolation. Observe
\begdi
	\dot{z}_1(t)=\sum_{i=0}^1 \bar{u}_i(t)\frac{1}{t}F_{x_i^{-1}(\eta_{{\mathbf s}_a})\bar{c}}[u](t)
			 =:e_2^Tz(t)\bar{u}_{i_1}(t)\frac{1}{t}.
\enddi
So the first row of $N_{i_1}$ is $e_2^T$, where $x_{i_1}$ is the first letter of $\eta_{\mathbf{s}_a}$, and the first row of the other realization matrix contains all zeroes. Continuing in this way,
\begdi
	\dot{z}_{k}(t)=\sum_{i=0}^1 \bar{u}_i(t)\frac{1}{t}F_{\eta_{{\mathbf s}_a}^{-1}(\eta_{{\mathbf s}_a})\bar{c}}[u](t)
			=:e_{k+1}^Tz(t)\bar{u}_{i_k}(t)\frac{1}{t}.
\enddi
Since in general $x_{i_k}=x_1$, the $k$-th row of $N_1$ is $e_{k+1}^T$, and the $k$-th row of the $N_0$ contains all zeroes. So far, this is in agreement with the proposed structure of the realization. Next observe that
\begdi
	\dot{z}_{k+1}(t)	=\sum_{i=0}^1 \bar{u}_i(t)\frac{1}{t}F_{x_i^{-1}(\bar{c})}[u](t)
				=\sum_{i=0}^1 \bar{u}_i(t)\frac{1}{t}\underbrace{F_{x_i^{-1}(\eta_{\mathbf{s}_c})}[u](t)}_{=:z_{k+2}(t)} +
		\sum_{j=0}^1 \bar{u}_j(t)\frac{1}{t}\underbrace{F_{x_j^{-1}(\eta_{\mathbf{s}_b}\bar{c})}[u](t)}_{=:z_{k+3}(t)}.
\enddi
In this way, new states are created until finally the term $F_{\bar{c}}[u](t)=z_{k+1}(t)$ reappears as it must. This produces an entry $\theta^{\abs{\mathbf{s}_b}}$ in $N_1$ and preserves the proposed structures of $N_0$ and $N_1$. But note, as in Example~\ref{ex:realization-L212n3}, that the process can continue to create new states, and the state $z_{k+1}(t)$ could reappear if $\eta_{{\mathbf s}_{c}}$ is a power of $\eta_{{\mathbf s}_b}$, a possibility that has not been excluded. In addition, this realization could produce {\em copies} of the the first $k$ states if $\eta_{{\mathbf s}_c}$ contains $\eta_{{\mathbf s}_a}$ as a factor. These copies will still preserve the desired structure, but this possibility points out that in general the final realization constructed by this process may not be minimal. Finally, the canonical embedding \rref{eq:brockett-embedding}, which is always needed if $\abs{s_c}>0$, yields the final elements of the proposed structure.
\endpr

Clearly, when non-periodic components are present, giving a precise general form of the matrices $N_0$ and $N_1$ is not as simple as in the purely periodic case.


\subsection{The dendriform setting}
\label{ssect:dend}

Recall that MPLs satisfy shuffle product identities, which are derived from integration by parts for the iterated integrals in \eqref{eq:MPLint}. For instance,
$$
	\mathrm{Li}_{(2)}(t)\mathrm{Li}_{(2)}(t)=4\mathrm{Li}_{(3,1)}(t) + 2\mathrm{Li}_{(2,2)}(t).
$$
In slightly more abstract terms this can be formulated using the notion of a dendriform algebra. Indeed, for any $t_0 < t_1$, the space $C[t_0,t_1]$ is naturally endowed with such a structure consisting of two products:
\begin{subequations}
\label{dend-example}
\begin{align}
	f \succ g &:= I(f)g	\label{righDend}\\
	f \prec g &:= fI(g),	\label{leftDend}
\end{align}
\end{subequations}
where $I$ is the Riemann integral operator defined by $I(f)(t,t_0):=\int_{t_0}^t f(s)\,ds$, and which are easily seen to satisfy the axioms of a {\it{dendriform algebra}}
\begin{subequations}
\label{dend-alg}
\begin{align*}
f\succ (g\succ h)&=(f*g)\succ h\\
(f\succ g)\prec h&=f\succ (g\prec h)\\
(f\prec g)\prec h&=f\prec (g*h),
\end{align*}
\end{subequations}
where
$$
	f*g := f\succ g + f\prec g
$$
is an associative product.
The example \eqref{dend-example} above moreover verifies the extra commutativity property $f\succ g=g\prec f$, making it a Zinbiel algebra\footnote{The space of continuous maps on $[t_0,t_1]$ with values in the algebra $\Cal M_n(\mathbb R)$ is also a dendriform algebra, with $\prec $ and $\succ $ defined the same way. But it is Zinbiel only for $n=1$.}
$$
	(f\prec g)\prec h = f \prec ( g \prec h + h \prec g).
$$
This is another way of saying that Chen's iterated integrals define a shuffle product, which gives rise to the shuffle algebra of MPLs. For more details, including a link between general, i.e., not necessarily commutative, dendriform algebras and Fliess operators, the reader is referred to \cite{Duffaut-Gray-Ebrahimi-Fard_14,Duffaut-Gray-Ebrahimi-Fard_16, EM_09}.

In the following, the focus is on the commutative dendriform algebra $(C[t_0,t_1], \succ, \prec)$. The linear operator $R^\succ _g:C[t_0,t_1]\to C[t_0,t_1]$ is defined for $g \in C[t_0,t_1]$ by right multiplication using \eqref{righDend}
$$
	R^\succ _g(f):=f\succ g.
$$
Now add the distribution $\delta=\delta_{t_0}$ to the dendriform algebra $C[t_0,t_1]$. In view of the identity $I(\delta)=1$ on the interval $[t_0,t_1]$, it follows that $R^\succ _f(\delta)=\delta\succ f=f$ for any $f\in C[t_0,t_1]$. Consider next the specific functions $u_0(t)=1/t$ and $u_1(t)=1/(1-t)$ which appeared above (with $t_0=0$ and $t_1=1$ here), and the corresponding linear operators $R^\succ _{u_0}$ and $R^\succ _{u_1}$. The notation $u_0=\wt x_0$ and $u_1=\wt x_1$ is useful. For any word $w = x_0^{s_1-1}x_1 \cdots x_0^{s_l-1}x_1 \in x_0X^*x_1$, the linear operator $R^\succ _w$ is defined as the composition of the linear operators associated to its letters, namely,
\begin{equation*}
	R^\succ _w=(R^\succ _{\wt x_0})^{s_1-1}R^\succ _{\wt x_1}\cdots (R^\succ _{\wt x_0})^{s_l-1}R^\succ _{\wt x_1}
\end{equation*}
for $w=w_1\cdots w_{|\mathbf s|}=x_0^{s_1-1}x_1\cdots x_0^{s_l-1}x_1$. Using the shorthand notation $R^\succ _w=R^\succ _{\mathbf s}$ with $\mathbf s=(s_1,\ldots,s_l)$, the multiple polylogarithm $\mop{Li}_{\mathbf s}$ obviously satisfies
\begin{align}\label{li-dend}
	\frac{d}{dt}\mop{Li}_{\mathbf s} &=R_{\mathbf s}^\succ (\delta).
\end{align}
From \eqref{li-dend} it follows immediately that
\begdi
\label{dend-periodic}
	\frac{d}{dt}\Cal L_{\mathbf s}(t,\theta)=\sum_{k=0}^\infty \theta^{k|\mathbf s|}(R_{\mathbf s}^\succ )^k(\delta),
\enddi
which in turn yields
\begin{equation}
\label{dend-periodic-bis}
	\frac{d}{dt}\Cal L_s(t,\theta)=\delta+\theta^{|\mathbf s|}R_{\mathbf s}^\succ \Big( \frac{d}{dt}\Cal L_s(t,\theta)\Big).
\end{equation}
Equation \eqref{dend-periodic-bis} is a dendriform equation of degree $(|\mathbf s|,0)$ in the sense of \cite[Section 7]{EM_09}. The general form of the latter is
\begin{equation}
\label{gen-dend-eq}
	X = a_{00} + \sum_{q=1}^{|\mathbf s|}\theta^{q}\sum_{j=1}^{q} (\cdots (X \succ a_{q1})\succ  a_{q1} \cdots )\succ a_{qq}
\end{equation}
with $a_{00}:=\delta$, $a_{qj}=0$ for $q<|\mathbf s|$ and $a_{|\mathbf s|j}:=\wt w_j$, matching the notations of equation (46) in reference \cite{EM_09}. The general solution $X$ of \eqref{gen-dend-eq} is the first coefficient of a vector $Y$ of length $|\mathbf s|$ whose coefficients (discarding the first one) are given by $\theta^jR^\succ_{w_1\cdots w_j}(X)$ for $j=1,\ldots,|\mathbf s|-1$. This vector satisfies the following matrix dendriform equation of degree $(1,0)$:
\begin{equation}
\label{matrix-dend-eq}
	Y=(\delta,\underbrace{0,\ldots,0}_{|\mathbf s|-1})+\theta Y\succ N,
\end{equation}
where the matrix\footnote{The size of the matrix can be reduced from $1+|\mathbf s|(|\mathbf s|-1)/2$ to $|\mathbf s|$ by eliminating rows and columns of zeroes due to the particular form of \eqref{dend-periodic-bis} compared to equation (46) in \cite{EM_09}.} $N$ is given by:
\begdi
N=
\begin{bmatrix}
	0 	&	\wt w_{1}	&	0	&	0	&\cdots&	0\\
	0	&	0	&	\wt w_{2}	&	0	&\cdots&	0\\
	0	&	0	&	0	&	\wt w_{3}	&\cdots&	0\\
	\vdots&\vdots	&\ddots	&\ddots	&\vdots&\\
	0	&	0	&	0	&	0	&\cdots&	\wt w_{|\mathbf s|-1}\\
\wt w_{|\mathbf s|}&0&0&0&\cdots&0
\end{bmatrix}.
\enddi
First, observe that the $|\mathbf s|$-fold product $(\cdots (N \succ N) \succ \cdots ) \succ N$ yields a diagonal matrix with the entry $\frac{d}{dt}\mop{Li}_{\mathbf s}(t)$ in the position $(1,1)$. Second, matrix $N$ splits into $N = N_0u_0 + N_1u_1$ with $N_0,N_1$ as in \eqref{eq:bilinear-Fetas}. Equation \eqref{matrix-dend-eq} essentially corresponds to the integral equation deduced from \eqref{eq:general-SISO-bilinear-system} giving the state $z(t)$.\\

\noindent 
The case with non-periodic components can also be handled in this setting. Observe
\begdi
\frac d{dt}\Cal L_{\mathbf s_a\{\mathbf s_b\}\mathbf s_c}=R_{\mathbf s_a}^\succ \left(\frac d{dt}\Cal L_{\{\mathbf s_b\}\mathbf s_c}\right),
\enddi
and the term $X'=\frac d{dt}\Cal L_{\{\mathbf s_b\}\mathbf s_c}$ satifies the dendriform equation
\begin{equation}\label{dend-np}
X'=R_{\mathbf s_c}^\succ(\delta)+\theta^{|\mathbf s_b|}R_{\mathbf s_b}^\succ(X').
\end{equation}
Equation \eqref{dend-np} is again a dendriform equation of degree $(|\mathbf s_b|,0)$ with $a_{00}=R_{\mathbf s_c}^\succ(\delta)$, $a_{qj}=0$ for $q<|\mathbf s_b|$ and $a_{|\mathbf s_b|j}=w_j$ using the notation in \cite{EM_09}. The general solution $X'$ of \eqref{dend-np} is the first coefficient of a vector $Y'$ of length $|\mathbf s_b|$ whose coefficients (discarding the first one) are given by $\theta^jR_{w_1\cdots w_j}^\succ(X')$ for $j=1,\ldots,|\mathbf s_b|-1$. This vector satisfies the following matrix dendriform equation of degree $(1,0)$
\begdi
Y'=(R_{\mathbf s_c}^\succ(\delta),\underbrace{0,\ldots,0}_{|\mathbf s_b|-1})+\theta Y'\succ M',
\enddi
where the matrix $M'$ is given by:
\begdi
M'=
\begin{bmatrix}
0 &\wt w_1&0&0&\cdots&0\\
0&0&\wt w_2&0&\cdots&0\\
0&0&0&\wt w_3&\cdots&0\\
\vdots&\vdots&\ddots&\ddots&\vdots\\
0&0&0&0&\cdots&\wt w_{|\mathbf s_b|-1}\\
\wt w_{|\mathbf s_b|}&0&0&0&\cdots&0
\end{bmatrix}.
\enddi

One can ask the question whether the term $X=\frac d{dt}\Cal L_{\mathbf s_a\{\mathbf s_b\}\mathbf s_c}$ itself is a solution of a dendriform equation. In fact, a closer look reveals that the theory of linear dendriform equations presented in \cite{EM_09} has not been sufficiently developed to embrace this more complex setting. In the light of Theorem \ref{th:realization-GF-nonperiodic-multiple-polylogs}, it is clear that the results in \cite{EM_09} should be adapted in order to address this question. Such a step, however, is beyond the scope of this paper and will thus be postponed to another work. It is worth mentioning that the matrix $N$ needed in the linear dendriform equation
$$
	Y' = (0,\delta,0,0,0,0,0)+ \theta Y' \succ N
$$
to match the result from Example \ref{ex:realization-L212n3} has the form

\begdi
N=
\begin{bmatrix}
	0 	&0		&0		&0		&0		&0		&0\\
	0	&0		&\wt w_1	&0		&0		&0		&0\\
	0	&0		&0		&\wt w_2	&0		&0		&0\\
	0	&0		&0		&0		&\wt w_3	&0		&0\\
	0	&0		&0		&0		&0		&\wt w_4	&0\\
	0	&0		&0		&0		&\wt w_6	&0		&\wt w_5\\
	\wt w_7 	&0		&0		&0		&0		&0		&0\\
\end{bmatrix},
\enddi
which reflects the canonical embedding of Brockett. The first component of the vector $Y'$ contains the solution. As indicated earlier, a proper derivation of this result in the context of general dendriform algebras, i.e.,  extending the results in \cite{EM_09}, lies outside the scope of the present paper.


\section{Examples}
\label{sec:examples}

In this section, three examples of the method described above are given corresponding to the generating functions behind the identities
\rref{eq:zeta2n-equals-ratio}, \rref{eq:zeta4n-equals-4zeta31n}, and \rref{eq:Hoffman-conjecture}.
\begin{figure}[h]
\begin{center}
	\includegraphics*[scale=0.4]{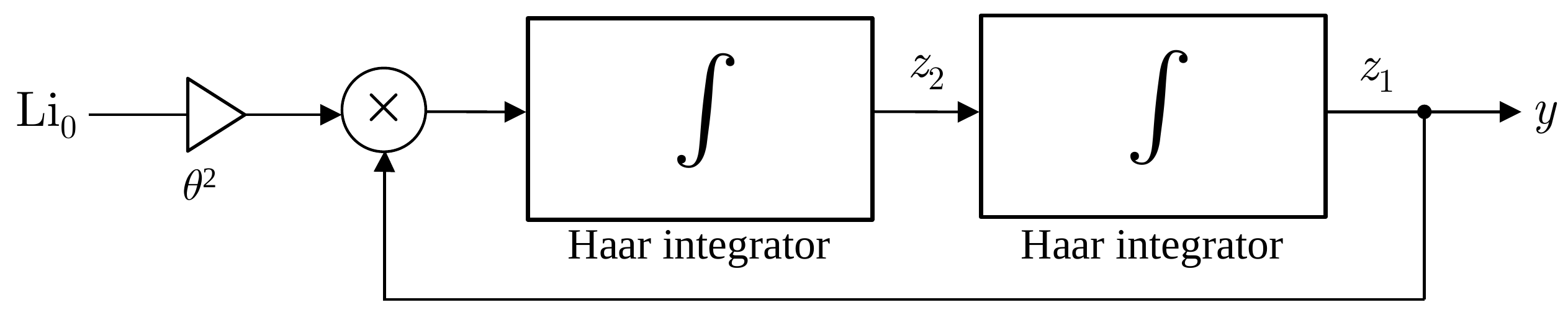}
\end{center}
\caption{Unity feedback system realizing $\mathcal{L}_{(2)}(t,1)$}
\label{fig:realization-generating-function-L2}
\end{figure}

\begin{exam}
Consider the generating function $\mathcal{L}_{(2)}(t,\theta)$. This example
is simple enough that a bilinear realization can be identified directly from \rref{eq:L2-ode}.
For any fixed $\theta$ define the first state variable to be $z_1(t)=\mathcal{L}_{(2)}(t,\theta)$,
and the second state variable to be $z_2(t)=t\,d\mathcal{L}_{(2)}(t,\theta)/dt$. In which case,
\begin{subequations}
\label{eq:bilinear-system-example-L2}
\begin{align}
\dot{z}_1(t)&=z_2(t)\frac{1}{t},\;\;z_1(0)=1 \\
\dot{z}_2(t)&=\theta^2\, z_1(t) \frac{t}{1-t}\,\frac{1}{t},\;\;z_2(0)=0 \\
y(t)&=z_1(t).
\end{align}
\end{subequations}
Thereupon, system \rref{eq:bilinear-system-example-L2} assumes the form of a bilinear system as given by
\rref{eq:bilinear-Fetas}, where the inputs are set to be $\bar{u}_0(t)=1$ and $\bar{u}_1(t)=\mathrm{Li}_0(t)=t/1-t$,
i.e.,
\begdi
N_0=N_0(2)=
\left[\begin{array}{cc}
0 & 1 \\
0 & 0
\end{array}\right],\;\;
N_1=N_1(2)=
\left[\begin{array}{cc}
0 & 0 \\
\theta^2 & 0
\end{array}\right],\;\;
z(0)=C^T=
\left[\begin{array}{c}
1 \\
0
\end{array}\right]
\enddi%
\begin{figure}[t]
\begin{center}
	\includegraphics*[scale=0.45]{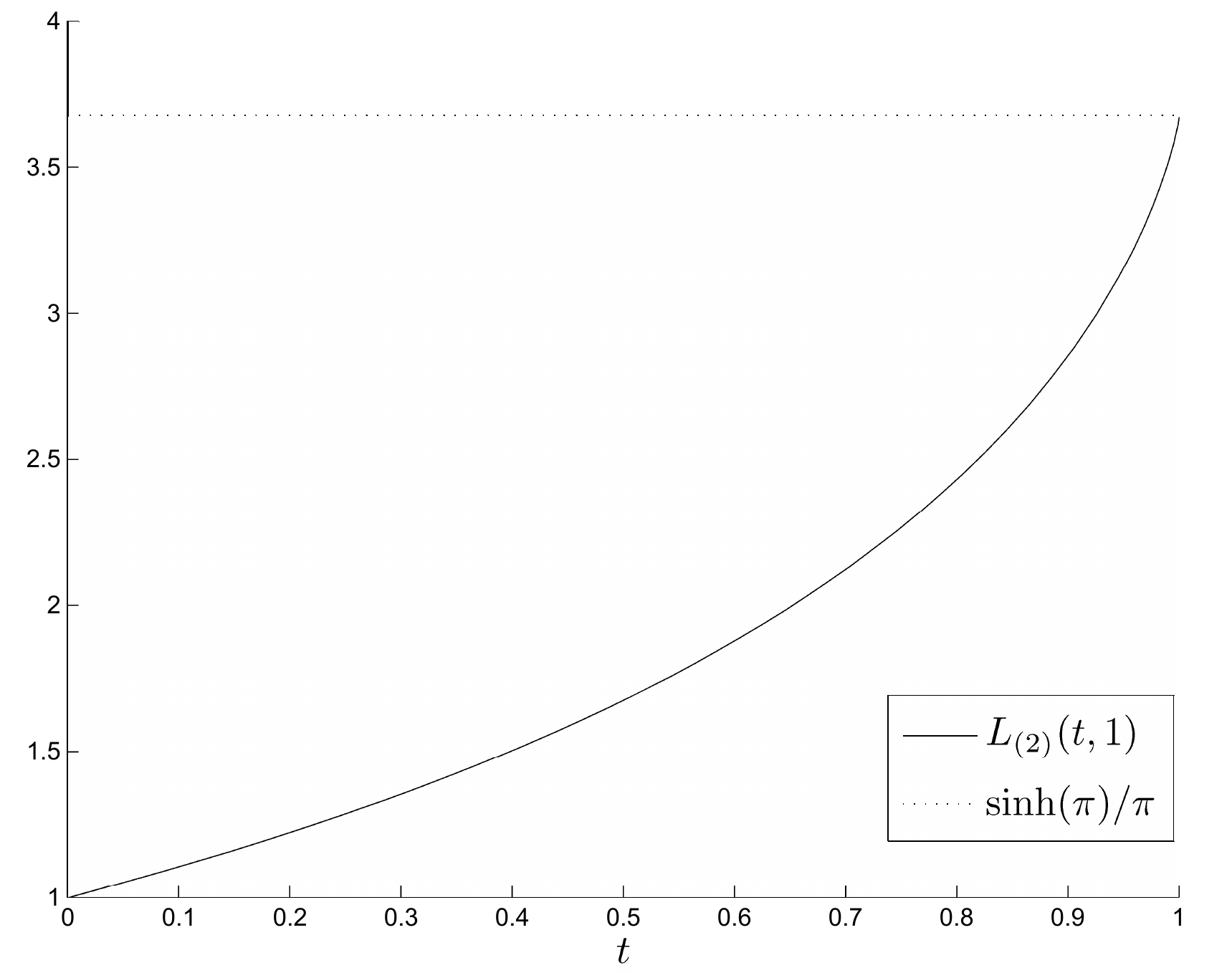}
\end{center}
\caption{Plot of $\mathcal{L}_{(2)}(t,1)$ versus $t$}
\label{fig:periodic-poylog-generating-function-Ls2-simulated}
\end{figure}%
(recall the $1/t$ factors in \rref{eq:bilinear-system-example-L2} are absorbed into Haar integrators).
A simulation diagram for this realization suitable for Matlab's Simulink simulation software is shown Figure~\ref{fig:realization-generating-function-L2}.
Setting $\theta=1$ and using Simulink's default integration routine {\tt ode45} (Dormand-Prince method \cite{Dormand-Prince_80}) with a variable step size lower
bounded by $10^{-8}$, Figure~\ref{fig:periodic-poylog-generating-function-Ls2-simulated} was generated showing
$\mathcal{L}_{(2)}(t,1)=F_{(x_0x_1)^\ast}[\mathrm{Li}_0](t)$
as a function of $t$.  In particular, it was found numerically that $\mathcal{L}_{(2)}(1,1)\approx 3.6695$,
which compares favorably to the theoretical
value derived from \rref{eq:zeta2n-equals-ratio}:
\begdi
\mathcal{L}_{(2)}(1,1)=\sum_{n=0}^\infty \zeta(\{2\}^n)=\sum_{n=0}^\infty \frac{\pi^{2n+1}}{(2n+1)^n}=\frac{\sinh(\pi)}{\pi}=3.6761.
\enddi
Better estimates can be found by more carefully addressing the singularities at the boundary conditions $t=0$ and $t=1$ in the Haar integrators.
\end{exam}

\begin{exam}
In order to validate \rref{eq:zeta4n-equals-4zeta31n}, the identity \rref{eq:L4-L31} is checked numerically. Since the generating functions
$\mathcal{L}_{(4)}$ and $\mathcal{L}_{(3,1)}$ are periodic, Theorem~\ref{th:realization-GF-periodic-multiple-polylogs} applies.
For $\mathbf{s}=(4)$ the corresponding bilinear realization is
\begdi
N_0=N_0(4)=
\left[
\begin{array}{cccc}
0 & 1 & 0 & 0 \\
0 & 0 & 1 & 0 \\
0 & 0 & 0 & 1 \\
0 & 0 & 0 & 0
\end{array}
\right],\;\;
N_1=N_1(4)=
\left[
\begin{array}{cccc}
0 & 0 & 0 & 0 \\
0 & 0 & 0 & 0 \\
0 & 0 & 0 & 0 \\
\theta^4 & 0 & 0 & 0
\end{array}
\right],\;\;
z(0)=C^T=
\left[\begin{array}{c}
1 \\
0 \\
0 \\
0
\end{array}\right].
\enddi
\begin{figure}[t]
\begin{center}
	\includegraphics*[scale=0.45]{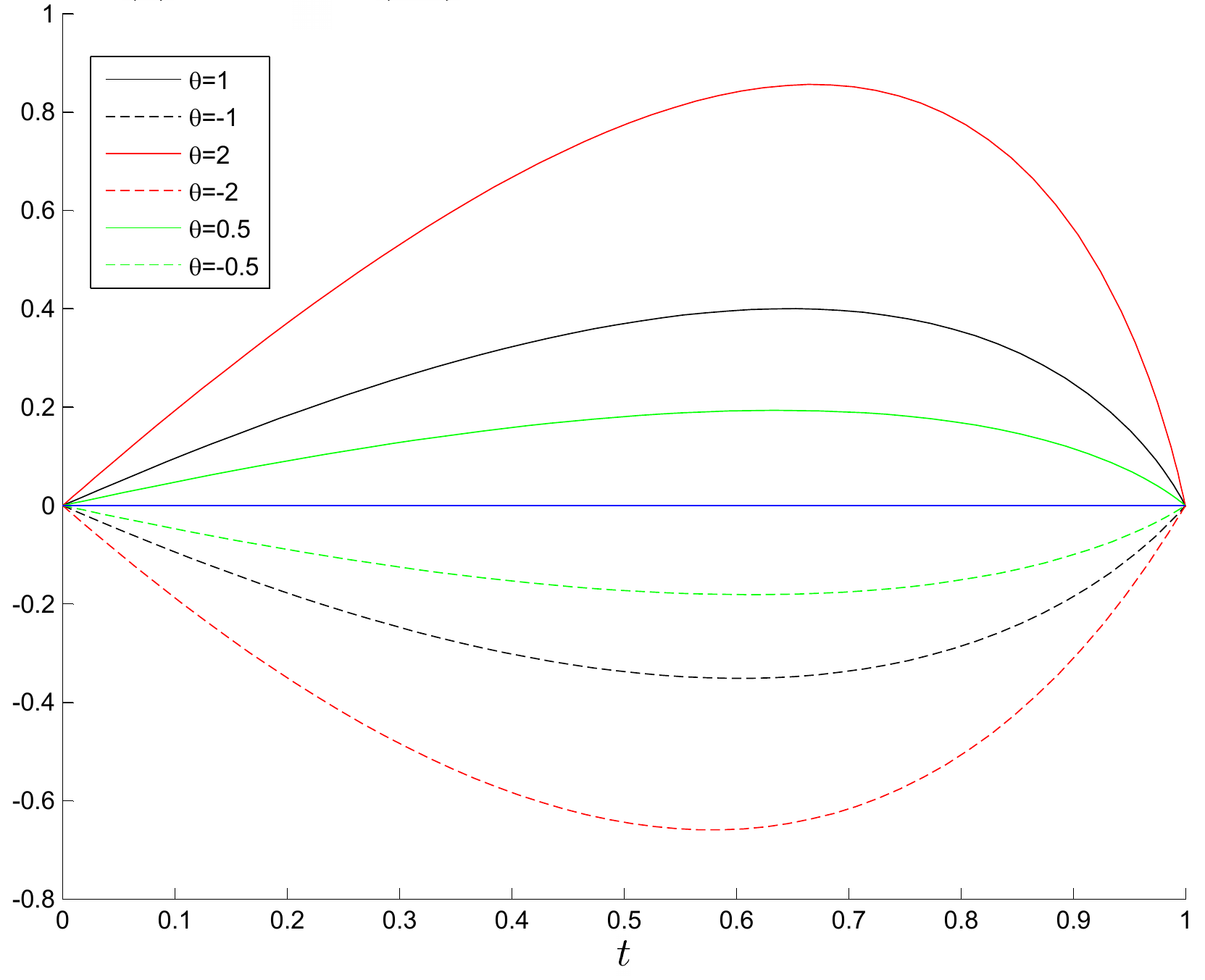}
\end{center}
\caption{Plot of $\mathcal{L}_{(4)}(t,\theta)-\mathcal{L}_{(3,1)}(t,\sqrt{2}\theta)$ versus $t$ for different values of $\theta$}
\label{fig:L4-L31-simulated}
\end{figure}
For $\mathbf{s}=(3,1)$ the bilinear realization is
\begdi
N_0=
\left[
\begin{array}{ccc|c}
0 & 1 & 0 & 0 \\
0 & 0 & 1 & 0 \\
0 & 0 & 0 & 0 \\ \hline
0 & 0 & 0 & 0
\end{array}
\right],\;\;
N_1=
\left[
\begin{array}{c|ccc}
0 & 0 & 0 & 0 \\
0 & 0 & 0 & 0 \\
0 & 0 & 0 & 1 \\ \hline
\theta^4 & 0 & 0 & 0
\end{array}
\right],\;\;
z(0)=C^T=
\left[\begin{array}{c}
1 \\
0 \\
0 \\
0
\end{array}\right].
\enddi
These two dynamical systems were simulated using Haar integrators in Simulink and the difference \rref{eq:L4-L31} was computed as a function of $t$ as shown in
Figure~\ref{fig:L4-L31-simulated}.
As expected, this difference is very close to zero
when $t=1$ no matter how the parameter $\theta$ is selected. This is pretty convincing numerical
evidence supporting \rref{eq:zeta4n-equals-4zeta31n}, which as discussed
in the introduction is known to be true.
\end{exam}

\begin{exam}
Now the method is applied to the generating functions behind the Hoffman conjecture \rref{eq:Hoffman-conjecture}.
In this case, each multiple polylogarithm has non-periodic components, so Theorem~\ref{th:realization-GF-nonperiodic-multiple-polylogs} has to be applied three times.
The realization for $\mathcal{L}_{(2,1,\{2\},3)}(t,\theta)$ was presented in Example~\ref{ex:realization-L212n3}. Following a similar approach, the realization for $\mathcal{L}_{(\{2\},2,2,2)}(t,\theta)$ and $\mathcal{L}_{(\{2\},3,3)}(t,\theta)$ are, respectively,
\begdi
N_0=
\left[
\begin{array}{cc|cccc|c}
0 & 1 & 0 & 0 & 0 & 0 & 0 \\
0 & 0 & 0 & 0 & 0 & 0 & 0 \\ \hline
0 & 0 & 0 & 1 & 0 & 0 & 0 \\
0 & 0 & 0 & 0 & 0 & 0 & 0 \\
0 & 0 & 0 & 0 & 0 & 1 & 0 \\
0 & 0 & 0 & 0 & 0 & 0 & 0 \\ \hline
0 & 0 & 0 & 0 & 0 & 0 & 0
\end{array}
\right],\;
N_1=
\left[
\begin{array}{cc|cccc|c}
0 & 0 & 0 & 0 & 0 & 0 & 0 \\
\theta^2 & 0 & 1 & 0 & 0 & 0 & 0 \\ \hline
0 & 0 & 0 & 0 & 0 & 0 & 0 \\
0 & 0 & 0 & 0 & 1 & 0 & 0 \\
0 & 0 & 0 & 0 & 0 & 0 & 0 \\
0 & 0 & 0 & 0 & 0 & 0 & 1 \\ \hline
0 & 0 & 0 & 0 & 0 & 0 & 0
\end{array}
\right],\;
z(0)=
\left[\begin{array}{c}
1 \\
0 \\ \hline
0 \\
0 \\
0 \\
0 \\ \hline
1
\end{array}\right],\;
C^T=
\left[\begin{array}{c}
1 \\
0 \\ \hline
0 \\
0 \\
0 \\
0 \\ \hline
0
\end{array}\right]
\enddi
and
\begdi
N_0=
\left[
\begin{array}{cc|cccc|c}
0 & 1 & 0 & 0 & 0 & 0 & 0 \\
0 & 0 & 1 & 0 & 0 & 0 & 0 \\ \hline
0 & 0 & 0 & 0 & 0 & 0 & 0 \\
0 & 0 & 0 & 0 & 1 & 0 & 0 \\
0 & 0 & 0 & 0 & 0 & 1 & 0 \\
0 & 0 & 0 & 0 & 0 & 0 & 0 \\ \hline
0 & 0 & 0 & 0 & 0 & 0 & 0
\end{array}
\right],\;
N_1=
\left[
\begin{array}{cc|cccc|c}
0 & 0 & 0 & 0 & 0 & 0 & 0\\
\theta^2 & 0 & 0 & 0 & 0 & 0 & 0 \\ \hline
0 & 0 & 0 & 1 & 0 & 0 & 0\\
0 & 0 & 0 & 0 & 0 & 0 & 0\\
0 & 0 & 0 & 0 & 0 & 0 & 0\\
0 & 0 & 0 & 0 & 0 & 0 & 1\\ \hline
0 & 0 & 0 & 0 & 0 & 0 & 0
\end{array}
\right],\;
z(0)=
\left[\begin{array}{c}
1 \\
0 \\ \hline
0 \\
0 \\
0 \\
0 \\ \hline
1
\end{array}\right],\;
C^T=
\left[\begin{array}{c}
1 \\
0 \\ \hline
0 \\
0 \\
0 \\
0 \\ \hline
0
\end{array}\right].
\enddi

\begin{figure}[t]
\begin{center}
	\includegraphics*[scale=0.45]{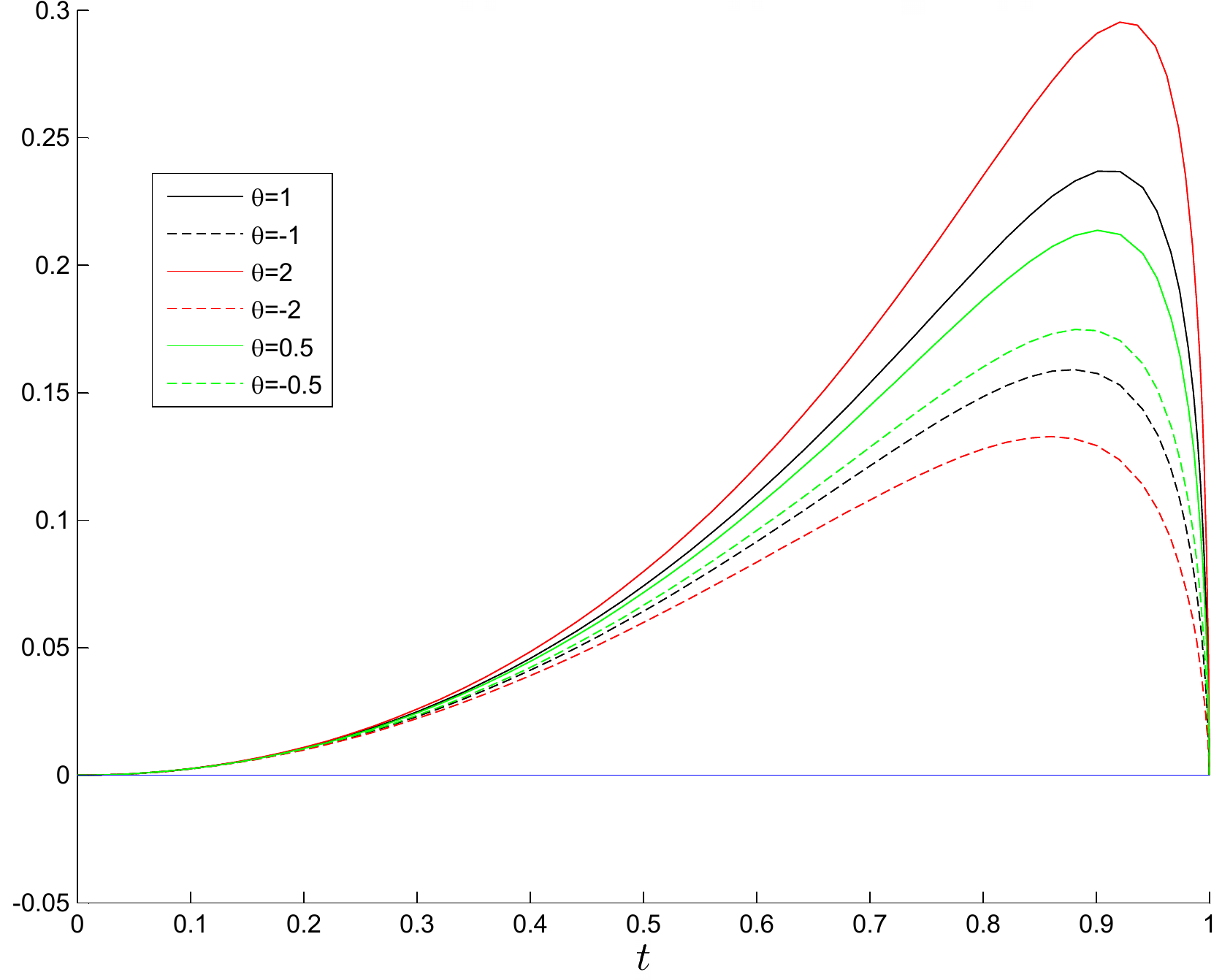}
\end{center}
\caption{Plot of $\mathcal{L}_{(\{2\},2,2,2)}(t,\theta)+2\mathcal{L}_{(\{2\},3,3)}(t,\theta)-\mathcal{L}_{(2,1,\{2\},3)}(t,\theta)$ versus $t$ for different values of $\theta$}.
\label{fig:Hoffman-MZV-conjecture-simulated}
\end{figure}%
These dynamical systems were simulated to estimate numerically the left-hand side of \rref{eq:GF-identity-Hoffman-conjecture} as shown in Figure~\ref{fig:Hoffman-MZV-conjecture-simulated}. As in the previous example, the case where $t=1$ is of primary interest. This value is again very close to zero for every choice of $\theta$ tested. It is highly likely therefore that the Hoffman conjecture is true.
\end{exam}


\section{Conclusions}
\label{sect:conclusion}

A systematic way was given to numerically evaluate the generating function of
periodic multiple polylogarithm using Chen--Fliess series with rational generating series. The method involved
mapping the corresponding Chen--Fliess series to a bilinear dynamical system, which could then be simulated numerically
using Haar integration.  A standard form for such a
realization was given, and the method was generalized to the case where the multiple polylogarithm could have
non-periodic components. The method was also described in the setting of dendriform algebras.
Finally, the technique was used to numerically validate the Hoffman conjecture.


\section*{Acknowledgements}

The first author is supported by Ram\'on y Cajal research grant RYC-2010-06995 from the Spanish government. The second author was supported by grant SEV-2011-0087 from the Severo Ochoa Excellence Program at the Instituto de Ciencias Matem\'{a}ticas in Madrid, Spain. This research was also supported by a grant from the BBVA Foundation.


\end{document}